\documentclass[english,11pt]{article}

\usepackage[latin1]{inputenc}
\usepackage{amsmath,amsthm,amssymb}
\usepackage{graphicx}
\usepackage{color}

\def\proof{{\it Proof.}\ }

\def\eq#1{(\ref{#1})}

\def\neweq#1{\begin{equation}\label{#1}}
\def\endeq{\end{equation}}

\def\ep{\varepsilon}

\def\RR{{\mathbb R} }
\def\NN{{\mathbb N} }

\def\ZZ{{\mathbb Z} }

\def\di{\displaystyle}
\def\ri{\rightarrow}

\def\kplus{{\cal K}_1(H)}
\def\kdoi{{\cal K}_2(H)}
\def\kdoimu{{\cal K}_{2,\mu}(H)}
\newtheorem{theorem}{Theorem}[section]

\title{\sc An infinite dimensional version of the Schur convexity property and applications}
\author{Claude VALL\'EE$^a$ and
Vicen\c tiu R\u ADULESCU$^b$\\
\small $^a$ Laboratoire de M\'ecanique des Solides, UMR 6610,
Universit\'e de Poitiers,\\
\small    Boulevard Marie et Pierre Curie,  BP 30179,
86962 Futuroscope Chasseneuil,  France\\ \small E-mail: {\tt
vallee@lms.univ-poitiers.fr}\\
\small $^b$ Department of Mathematics, University of Craiova,
200585 Craiova, Romania\\ \small E-mail: {\tt vicentiu.radulescu@math.cnrs.fr}}

\date{}

\begin{document}
\baselineskip16pt \maketitle
\renewcommand{\theequation}{\arabic{section}.\arabic{equation}}
\catcode`@=11 \@addtoreset{equation}{section} \catcode`@=12

\begin{abstract} We extend to infinite dimensional separable Hilbert spaces the Schur convexity 
property
of eigenvalues of a symmetric matrix with real entries. Our framework includes both the case of 
linear, selfadjoint, compact operators, and that of linear selfadjoint operators that can be 
approximated by operators of finite rank
and having a countable family of eigenvalues. The abstract
results of the present paper are illustrated  by several examples from mechanics or quantum 
mechanics, including the Sturm-Liouville problem, the Schr\"odinger equation, and the harmonic 
oscillator.
\end{abstract}

\section{Introduction}
An important notion in the finite dimensional theory of convex functions is that of {\it Schur 
convexity}. 
Roughly speaking, Schur-convex functions are real-valued mappings which
are monotone with respect to the majorization ordering. A rigorous definition is stated in what 
follows.
Let $\RR^n_\geq$ denote the cone of vectors with nonincreasing components, that is,
$$\RR^n_\geq=\{x=(x_1,x_2,\ldots ,x_n)\in\RR^n;\ x_1\geq x_2\geq \ldots\geq x_n\}\,.$$
The dual cone of the cone $\RR^n_\geq$ is defined by
$$(\RR^n_\geq)^+=\left\{ y\in\RR^n;\  (x,y) \geq 0\mbox{ for all }x\in\RR^n_\geq \right\}\,.$$
A straightforward computation shows that
$$(\RR^n_\geq)^+=\left\{ y\in\RR^n;\  \sum_{i=1}^j y_i\geq 0\mbox{ for all }j=1,\ldots ,n-1 \mbox{ 
and }
 \sum_{i=1}^ny_i=0 \right\}\,.$$

We recall (see, e.g., Roberts and Varberg \cite{roberts}, Borwein and Lewis \cite{borlew}) that a 
function $f:\RR^n\ri\RR$ is {\it Schur convex}
if it is $(\RR^n_\geq)^+$-isotone, that is,
$$x,y\in \RR^n_\geq,\ y-x\in (\RR^n_\geq)^+\Longrightarrow f(x)\leq f(y)\,.$$

The Schur-convex functions were introduced by Schur \cite{schur23} in 1923  and they have many
important applications in analytic inequalities. Hardy, Littlewood and P\'olya \cite{hardylp}
were also interested in some
inequalities that are related to Schur-convex functions.
The notion of Schur-convexity has shown its importance in many domains. For instance, Merkle 
proved in \cite{merkle}
that if $I\subset\RR$ is an interval and $f:I\ri\RR$ is differentiable, then $f'$ is convex
if and only if the mapping
$$F(x,y)=\left\{
\begin{array}{lll}
&\di \frac{f(y)-f(x)}{y-x}\qquad &\di\mbox{if $y\not=x$}\\
& f'(x)\qquad &\di\mbox{if $y=x$}
\end{array}\right.
$$
is Schur convex. This property is applied in order to obtain some
inequalities for the ratio of Gamma functions. We also refer to
Hwang and  Rothblum \cite{hwa}, who study
optimization problems over partitions of a finite set and
obtain conditions that allow for simple constructions of partitions that are uniformly optimal for 
all Schur convex functions. Stochastic Schur convexity properties have been established by Shaked, 
Shanthikumar and  Tong 
\cite{shaked}. Exciting results such as
Schur's analytic criteria for Schur convexity, equivalence with Muirhead's inequality, 
majorization and stochastic matrix conditions in $\RR^n$, and Schur's majorization inequality
can be found in the excellent book by Steele \cite{steele}. Recently, Guan \cite{guan} has proved 
that the complete
elementary symmetric function $c_r=c_r(x)=\sum_{i_1+\ldots +i_n=r}x_1^{i_1}\ldots x_n^{i_n}$ and 
the function
$c_r(x)/c_{r-1}(x)$ are Schur-convex functions in $\RR^n_+=\{(x_1,\ldots ,x_n);\ x_i>0\}$, where 
$r$ is a positive 
integer and $i_1,\ldots ,i_n$ are nonnegative integers.

Zhang \cite{zhangx} proved that every Schur-convex function $f:D\subset\RR^n\ri\RR$ is a symmetric 
function,
that is, $f(x)=f\left(x_{\sigma (1)},\ldots ,x_{\sigma (n)}\right)$ 
for any permutation $\sigma\in{\cal P}_n$ and for all $x=(x_1,\ldots ,x_n)\in D$. 
The converse is not true (see, e.g.,
Roberts and Varberg \cite[p.~258]{roberts}). However, if $I$ is an open interval and $f:I^n\ri\RR$ 
is symmetric and 
of class $C^1$, then  $f$ is Schur-convex if and only if
$$(x_i-x_j)\left( \frac{\partial f}{\partial x_i}-\frac{\partial f}{\partial x_j}\right)\geq 
0\qquad\mbox{on $I^n$},$$
for all $i,j\in\{1,\ldots ,n\}$ (see Roberts and Varberg \cite[p.~259]{roberts}).

Eigenvalues of real
symmetric matrices exhibit remarkable convexity properties.
Let ${\bf S}^n$ denote the set of all 
symmetric matrices $X\in {\cal M}_{n,n}(\RR)$. 
In Borwein and Lewis \cite[p.~108]{borlew} it is stated the following elementary property of 
eigenvalues of $X\in {\bf S}^n$.

\smallskip
{\bf The Schur Convexity Property.} {\sl Let
$\lambda_1(X)\geq \lambda_2(X)
\geq\ldots\geq \lambda_n(X)$ be the eigenvalues (counted by multiplicity)
of an arbitrary matrix $X\in {\bf S}^n$. Assume that $\mu =(\mu_1,\mu_2,\ldots ,\mu_n)\in 
\RR^n_\geq$. 
Then the functional $\varphi (X)=\sum_{i=1}^n\mu_i\lambda_i(X)$ is sublinear.
}

\smallskip
A direct consequence of this result is that the mapping $\varphi :{\bf S}^n\ri\RR$ is convex.

In the particular case $\mu_1=\ldots =\mu_k=1$, $\mu_{k+1}=\ldots =\mu_n=0$ 
($1\leq k\leq n$), we deduce that the sum of the largest $k$  eigenvalues of a matrix $X\in {\bf 
S}^n$ is a convex function. 
An alternative proof is based on the observation that, for any fixed $1\leq k\leq n$,
\begin{equation}\label{01}
\lambda_1(X)+\lambda_2(X)+\ldots +\lambda_k(X)=\sup_{A\in{\cal A}}{\rm tr}\, (AA^TX)\,,
\end{equation}
where $${\cal A}=\left\{A\in{\cal M}_{n,k}(\RR);\ A^TA={\rm I}_k\right\}\,.$$
Since ${\cal A}$ is a compact set, the supremum in \eq{01} is attained in ${\cal A}$.
We deduce that the mapping ${\bf S}^n\ni X\longmapsto\lambda_1(X)+\lambda_2(X)+\ldots 
+\lambda_k(X)$
is convex, as a supremum of linear functions on ${\bf S}^n$.
The extreme situations $k=1$ and $k=n$ show that both the
largest eigenvalue of $X$ and the trace of $X$ are convex functions on ${\bf S}^n$. We also 
deduce, by taking differences, that $\sum_{j=k+1}^n\lambda_j(X)$ is a concave function, for all 
$1\leq k\leq n-1$. In particular, the mapping ${\bf S}^n\ni X\longmapsto\lambda_n(X)$ is concave. 

A classical result (see, e.g., Borwein and Lewis \cite{borlew}, and
Rockafellar and  Wets \cite{rock}) asserts that Schur convex functions are precisely
restrictions to $\RR^n_\geq$ of symmetric convex functions.
This result is strictly related to the class of convex functions $f:{\bf S}^n\ri\RR$
(like the functions $\sum_{j=1}^k\lambda_j(X)$) depending only on the eigenvalues of $X$.  
In fact, if we write diag$\,(\lambda)$
(where $\lambda=(\lambda_1,\ldots ,\lambda_n)\in\RR^n$) 
for the diagonal matrix with diagonal entries $\lambda_1,\ldots ,\lambda_n$, and define a function 
$\Phi:\RR^n\ri\RR$
by $\Phi (\lambda)=f({\rm diag}\,(\lambda))$, then $\Phi$ is convex and symmetric: $\Phi 
(\lambda)=
\Phi (\sigma\circ\lambda)$ for all permutation $\sigma\in{\cal P}_n$. The converse is also true: 
if $\Phi:\RR^n\ri\RR$
is a symmetric convex function then the function $f:{\bf S}^n\ri\RR$ defined by $f(X)=\Phi 
(\lambda(X))$ 
(where $\lambda (X)=(\lambda_1(X),\ldots ,\lambda_n(X))^T$) is convex and satisfies
$f(U^*XU)=f(X)$ whenever $U\in {\cal M}_{n,n}(\RR)$ is a unitary matrix. The above result is due 
to
Davis \cite{davis}. 

The above considerations show that it is natural to impose an adequate ``symmetry" assumption in 
order to obtain a
Schur convexity property for linear operators defined on arbitrary Hilbert spaces. That is why we 
consider throughout this paper linear selfadjoint operators defined on infinite dimensional 
Hilbert spaces.

\section{A Schur convexity property in Hilbert spaces}
In the first part of this section we establish an infinite dimensional version of the Schur 
convexity property
for linear, selfadjoint and compact operators defined on separable Hilbert spaces. Next, we extend 
this 
property to the class of linear selfadjoint operators that can be approximated by operators of 
finite rank. Several
examples from mechanics and quantum mechanics illustrate both cases.

\subsection{Schur convexity property for selfadjoint, compact operators}
Let $H$ be a separable Hilbert space and assume that $S:H\ri H$ is a linear, selfadjoint and 
compact operator. 
Since $S$ is compact then, by the Riesz-Schauder theorem (Theorem VI.15 in Reed and Simon 
\cite{reed}), the spectrum $\sigma (S)$
of $S$ is a discrete set having no limit points except perhaps the origin. Moreover, any 
$\lambda\in\sigma (S)\setminus
\{0\}$ is an eigenvalue of finite multiplicity. Next, the classical spectral theory of compact 
selfadjoint operators
(see, e.g., Brezis \cite[Proposition VI.9]{brezis}) ensures that $\sigma (S)\subset [m,M]$ and $m, 
M\in\sigma (S)$,
where $m=\inf\{ (Su,u);\ u\in H,\ \|u\|=1\}$ and $M=\sup\{ (Su,u);\ u\in H,\ \|u\|=1\}$. In 
conclusion, the spectrum of $S$ is discrete and it consists of a countable family of eigenvalues 
$\left(\lambda_n(S)\right)_{n\geq 1}$ with the additional
property that $\lambda_n(S)\ri 0$ as $n\ri\infty$. At this stage,
the Hilbert-Schmidt theorem (Theorem VI.16 in Reed and Simon \cite{reed}) 
implies that there is a complete orthonormal
basis $(e_n)_{n\geq 1}$ of $H$ such that $Se_n=\lambda_ne_n$ for all $n\geq 1$, where 
$\lambda_n=\lambda_n(S)$. So,
$Sx=\sum_{n=1}^\infty \lambda_n(x,e_n)e_n$, for all $x\in H$.

We observe that for any fixed positive integer $n$, the set 
$$\left\{\lambda\in\sigma (S);\ |\lambda|\geq\frac 1n\right\}$$
is either empty or finite. Thus, we can rearrange the eigenvalues of $S$ such that 
\begin{equation}\label{laspectr}
\lambda_1(S)\geq \lambda_2(S)\geq\ldots\geq \lambda_n(S)\geq\ldots >0>\ldots 
\geq\lambda_{-n}(S)\geq\ldots
\geq\lambda_{-2}(S)\geq\lambda_{-1}(S)\,.\end{equation}
Moreover, the unique limit point of the sequence $\left(\lambda_n(S)\right)_{n\in\ZZ}$ is 0.  
If $S$ has a finite number of negative eigenvalues (say, $n$), we denote them by
$\lambda_{-1}(S)\leq \ldots\leq \lambda_{-n}(S)$ and we set $\lambda_{-k}(S)$ for all $k\leq n+1$. 
We make a similar convention if $S$ has finitely many positive eigenvalues.
If $0$ is an eigenvalue of $S$, we denote $\lambda_0 (S)=0$.

Denote by $\kplus$ the vector space of linear, selfadjoint and compact operators $S:H\ri H$.

We prove the following infinite dimensional variant of  Schur's convexity property.

\begin{theorem}\label{th1}
Let $H$ be a separable Hilbert space and assume that $S:H\ri H$ is an arbitrary compact 
selfadjoint operator. 
Assume that the eigenvalues of $S$ are arranged as in \eq{laspectr} and let 
$(\mu_n)_{n\in \ZZ}$ be real numbers such that $\mu_1\geq\mu_2\geq\ldots\geq\mu_n\geq\ldots\geq 
\mu_{-n}\geq\ldots
\geq\mu_{-2}\geq\mu_{-1}$ and $\sum_{n=-\infty}^\infty
\mu_n$ is an absolutely convergent series.

Then the functional $\psi :\kplus\ri\RR$ defined by $\psi 
(S)=\sum_{n=-\infty}^\infty\mu_n\lambda_n(S)$ is convex
and lower semicontinuous.
\end{theorem}

\proof 
We first observe that since $S\in\kplus$ is not assumed to be a nuclear operator, then the series 
$\sum_{n\in\ZZ}
\lambda_n(S)$ is not necessarily convergent. However, our hypothesis that the series 
$\sum_{n=-\infty}^\infty
\mu_n$ is absolutely convergent implies that the series $\sum_{n=-\infty}^\infty\mu_n\lambda_n(S)$ 
is absolutely convergent, too,
so the mapping $\psi$ is well-defined. Indeed, for all $S\in\kplus$,
$$|\psi (S)|\leq \sum_{n=-\infty}^\infty |\mu_n|\cdot |\lambda_n(S)|\leq
\max\left\{-\lambda_{-1}(S),\lambda_1(S)\right\}\sum_{n=-\infty}^\infty |\mu_n|<\infty\,.$$

Any operator $S\in\kplus$ is the norm limit of a sequence of operators of finite rank. Indeed, if 
$(e_n)_{n\in\ZZ}$
is a complete orthonormal basis of $H$ so that $Se_n=\lambda_n(S)e_n$ for all $n\in\ZZ$, with 
$\lambda_n(S)$
arranged as in \eq{laspectr}, then 
$Sx=\sum_{n=-\infty}^\infty \lambda_n(S)(x,e_n)e_n$, for all $x\in H$. Set, for any $m\geq 1$,
$S_mx=\sum_{j=-m}^m \lambda_j(S)(x,e_j)e_j$, for all $x\in H$. Then $S_m\ri S$ in $L(H)$ as 
$m\ri\infty$ and the (nontrivial) eigenvalues of $S_m$ are 
$\lambda_1(S)\geq \ldots\geq \lambda_m(S)>0>\lambda_{-m}(S)\geq\ldots\geq\lambda_{-1}(S)$. 
So, by the finite-dimensional Schur convexity property,
the mapping
$$\psi_m :\kplus\ri\RR\,,\qquad\qquad \psi_m (S)=\sum_{j=-m}^m\mu_j\lambda_j(S)$$ is sublinear. 
So, for any $S,\, T\in\kplus$ and all $\alpha\in\RR$,
\begin{equation}\label{1}
\psi_m(S+T)\leq\psi_m(S)+\psi_m(T)\qquad\mbox{and}\qquad\psi_m(\alpha S)=|\alpha|\,\psi_m(S)\,.
\end{equation}

On the other hand,
$$\left| \psi (S)-\psi_m(S)\right|=\left| \sum_{|j|\geq m+1}
\mu_j\lambda_j(S) \right|\leq
\max\left\{-\lambda_{-m-1}(S),\lambda_{m+1}(S)\right\}\sum_{|j|\geq m+1} |\mu_j|\,.$$
Therefore 
\begin{equation}\label{2}\psi_m(S)\ri\psi (S)\qquad \mbox{as $m\ri\infty$}.\end{equation}
Thus, by \eq{1} and \eq{2}, $\psi$ is a sublinear functional. In particular, $\psi$ is convex.

It remains to argue that $\psi$ is lower semicontinuous, that is, $\psi 
(S)\leq\liminf_{n\ri\infty}
\psi (S_n)$ for all $S\in\kplus$, provided $S_n\in\kplus$ and $\|S_n-S\|\ri 0$ as $n\ri\infty$.
The key ingredient is Theorem 4.2 in Gohberg-Krein \cite{gohkre}, which asserts that
$\lambda_j(S)=\lim_{n\ri\infty}\lambda_j(S_n)$. Fix an integer $m\geq 1$ and choose arbitrarily
$0<\ep <\max\left\{-\lambda_{-m}(S),\lambda_m(S)\right\}$. 
It follows that there exists $N_0=N_0(\ep)\in\NN$ such that,
for all $n\geq N_0$,
$$\begin{array}{ll}\di \psi_m(S)&\di=\sum_{j=-m}^m\mu_j\lambda_j(S)=
\sum_{j=-m}^m\mu_j^+\lambda_j(S)-\sum_{j=-m}^m\mu_j^-\lambda_j(S)\\
&\di\leq \sum_{j=-m}^m\mu_j^+\left(\lambda_j(S_n)+\ep\right)-
\sum_{j=-m}^m\mu_j^-\left(\lambda_j(S_n)-\ep\right)\\ &\di 
=\sum_{j=-m}^m\mu_j\lambda_j(S_n)+\ep\sum_{j=-m}^m|\mu_j|=
\psi_m(S_n)+\ep \sum_{j=-m}^m|\mu_j|\,.\end{array}$$
Taking $\ep\ri 0$ we obtain $\psi_m(S)\leq \psi_m(S_n)$, for all positive integers $m$ and $n$.
So, for all $n\geq 1$,
$$\psi (S)=\lim_{m\ri\infty} \psi_m(S_n)\leq\lim_{m\ri\infty} \psi_m(S_n)=\psi (S_n)\,.$$
We deduce that $\psi (S)\leq\liminf_{n\ri\infty}
\psi (S_n)$ and the proof is concluded. \qed

\medskip
{\bf Examples.} 1. {\it Sturm-Liouville differential operators}. Many eigenvalue problems in 
quantum mechanics as well as classical physics are described by the Sturm-Liouville problem
\begin{equation}\label{sturm}
\left\{
\begin{array}{ll}
&\di -\frac{d}{dx}\left(p(x)\frac{dy}{dx}\right)+q(x)y=\Lambda y\qquad \mbox{in $(0,L)$}\\
&\di y(0)=y(L)=0\,,
\end{array}\right.
\end{equation}
where $y(x)$ is the quantum mechanical wave function or other physical quantity, while $p\in 
C^1[0,L]$
($p>0$ in $[0,L]$) and $q\in C[0,L]$ are given functions that are determined by the nature of the 
system of interest.
We can assume, without loss of generality, that $q\geq 0$ in $[0,L]$. Indeed, if not, we choose 
$C\in\RR$
sufficiently large such that $q+C\geq 0$ in $[0,L]$ (in such a case, $\Lambda$ is replaced by 
$\Lambda +C$ in
\eq{sturm}). Fix $f\in L^2(0,L)$. Thus, by the Lax-Milgram lemma, there exists a unique
$u\in H^2(0,L)\cap H^1_0(0,L)$ such that
$$
\left\{
\begin{array}{ll}
&\di -\frac{d}{dx}\left(p(x)\frac{dy}{dx}\right)+q(x)y=f\qquad \mbox{in $(0,L)$}\\
&\di y(0)=y(L)=0\,.
\end{array}\right.
$$
Let $S:L^2(0,L)\ri L^2(0,L)$ be the operator defined by $Sf=u$. Then, by Theorem VIII.20 in Brezis 
\cite{brezis}, 
$S$ is linear, selfadjoint, compact, and
nonnegative. Let $\lambda_1(S)\geq\lambda_2(S)\geq\ldots\geq \lambda_n(S)\geq\ldots >0$ denote the 
eigenvalues of
$S$. Then $\Lambda_n(S)=1/\lambda_n(S)$ is an eigenvalue corresponding to the Sturm-Liouville 
problem \eq{sturm}.
In the particular case $p\equiv 1$ and $q\equiv 0$, a straightforward computation shows that 
$\lambda_n(S)=L^2(n^2\pi^2)^{-1}$.

Let $\mu_n$ ($n\geq 1$) be real numbers such that $\mu_i\geq \mu_j$ if $i<j$
and such that the series $\sum_{n=1}^\infty\mu_n$ converges absolutely. So, by Theorem \ref{th1}, 
the mapping
$$S\longmapsto\sum_{n=1}^\infty
\mu_n\lambda_n(S)$$
is convex and lower semicontinuous.

\smallskip
2. {\it The electron atom model}. On the Hilbert space $H=L^2(\RR^3)$, let $x$, $y$, $z$ be the 
components of
the momentum of the electron and denote by $r=(x,y,z)$ its position. Consider on $H$ the 
selfadjoint operator
$$S=\Delta +\frac{\alpha}{|r|}\,,\qquad  \qquad |r|=\sqrt{x^2+y^2+z^2}\,.$$
Notice that the potential $V(|r|)=\alpha /|r|$ is the energy of the electric field surrounding the 
electron, $\alpha$
depends on the electron's charge, and $|r|$ is its distance from the atom's nucleus.
As established in Reed and Simon \cite{reed4}, $S$ has no eigenvalues for any $\alpha<0$ and, if 
$\alpha>0$,
then all eigenvalues of $S$ are
$$\lambda_n(S)=\frac{\alpha}{4n^2}\,,\qquad \qquad n=1,2,\ldots$$

Let $(\mu_n)_{n\geq 1}$ be a sequence of real numbers such that $\mu_1\geq 
\mu_2\geq\ldots\geq\mu_n\geq\ldots$
and the series $\sum_{n=1}^\infty\mu_n$ converges absolutely. So, by Theorem \ref{th1}, the 
mapping
$$S\longmapsto\sum_{n=1}^\infty
\mu_n\lambda_n(S)$$
is convex and lower semicontinuous. 

\smallskip
3. {\it Nonrelativistic model for $2$-electron atom}.
Set $H=L^2(\RR^6)$ and define on $H$ the selfadjoint operator
$$S=\Delta_1 +\frac{\alpha}{|r_1|}+\Delta_2 +\frac{\beta}{|r_2|}\,,$$
where $\alpha, \beta>0$, $r_k=(x_k,y_k,z_k)$,  
and 
$$\Delta_k=\frac{\partial^2}{\partial x_k^2}+\frac{\partial^2}{\partial y_k^2}+
\frac{\partial^2}{\partial z_k^2}\,,\qquad \mbox{for all $k=1,2$}.$$

Cf. Reed and Simon \cite{reed4}, the eigenvalues of $S$ are precisely
$$\lambda_{n,m}(S)=\frac{\alpha}{4n^2}+\frac{\beta}{4m^2}\,,\qquad\qquad n,m=1,2,\ldots$$

The countable family of positive numbers $(\lambda_{n,m}(S))_{n,m\geq 1}$ can be rearranged in a 
sequence
$(\gamma_p(S))_{p\geq 1}$ such that $\gamma_i(S)\geq\gamma_j(S)$ provided $i<j$.
Let $(\mu_p)_{p\geq 1}$ be a sequence of real numbers such that $\mu_1\geq 
\mu_2\geq\ldots\geq\mu_p\geq\ldots$
and the series $\sum_{p=1}^\infty\mu_p$ converges absolutely. 
Thus, by Theorem \ref{th1}, the mapping
$$S\longmapsto\sum_{p=1}^\infty
\mu_p\gamma_p(S)$$
is convex and lower semicontinuous.

\smallskip
4. {\it Schr\"odinger operators with periodic potential}. 
The basic equation of quantum mechanics is the Schr\"odinger equation
\begin{equation}\label{erwin1}
i\hbar\psi_t=-\frac{\hbar^2}{2m}\,\Delta\psi+V(x)\psi\,.\end{equation}
Schr\"odinger \cite{schre} studied the stationary equation
\begin{equation}\label{erwin2}
\lambda\varphi=-\frac{\hbar^2}{2m}\,\Delta\varphi+V(x)\varphi\,,\end{equation}
which follows from \eq{erwin1} through $\psi(x,t)=\varphi (x)e^{-i\lambda t/\hbar}$. From 
\eq{erwin2}
Schr\"odinger derived the spectrum of the hydrogen atom. In this case, $V$ is the potential of
the electrostatic attracting force of the atomic nucleus, while from the eigenvalues $\lambda$
of \eq{erwin2} one obtains the energy levels of the electron of the hydrogen atom. 

Solutions of Schr\"odinger's equation have to fulfill strict conditions to be useful in describing 
the electron. Some of the solutions are associated with special values of the electron's energy 
level, known as eigenvalues.
We consider in what follows the class of piecewise continuous potential
functions $V:\RR\ri\RR$ which are periodic of period $2\pi$. Let $S$ denote 
the one dimensional Schr\"odinger operator associated to $V$ defined on $L^2_{\rm per}(\RR)$ 
 with $2\pi$-periodic conditions. 
This operator is defined as follows: for any $f\in L^2_{\rm per}(\RR)$ periodic of period $2\pi$, 
let $u\in H^1_{\rm per}(\RR)$ be the unique solution of the problem 
$$
\left\{
\begin{array}{ll}
&\di -u''+V(x)u=f\qquad \mbox{in $(0,2\pi)$}\\
&\di u(0)=u(2\pi),\ u'(0)=u'(2\pi)\,.
\end{array}\right.
$$
Then $S$ is defined by $L^2_{\rm per}(\RR)\ni f\longmapsto u=Sf\in L^2_{\rm per}(\RR)$.
 According to Theorem XIII.89 in Reed and Simon \cite{reed4}, $S$ has a countable family of 
eigenvalues
$\lambda_1(S)>\lambda_2(S)>\ldots >\lambda_n(S)>\ldots$ and $\lambda_n(S)\ri 0$ as $n\ri\infty$.
Assume that $\mu_n$ ($n\geq 1$) are real numbers such that $\mu_i\geq \mu_j$ if $i<j$
and such that the series $\sum_{n=1}^\infty\mu_n$ converges absolutely. So, by Theorem \ref{th1}, 
the mapping
$$S\longmapsto\sum_{n=1}^\infty\mu_n\lambda_n(S)$$
is convex and lower semicontinuous.

\smallskip
5. {\it Indefinite weight elliptic problems on the whole space}. Consider the class of measurable 
functions
$V:\RR^N\ri\RR$ ($N\geq 3$) such that $V^+\in L^{N/2}(\RR^N)$, where $V=V^+-V^-$. We observe that 
this class contains
potentials $V$ satisfying $V^+(x)\leq C(1+|x|^2)^{-\alpha}$ for all $x\in\RR^N$, where $\alpha >1$ 
and $C$
is a positive constant. For some fixed $\lambda>0$, let $E$ be the completion of 
$C^\infty_0(\RR^N)$ with
respect to the norm
$$\|u\|^2=\int_{\RR^N}\left[|\nabla u|^2+\max \left(\lambda V^-,\omega\right)u^2\right]dx\,,$$
where $\omega (x)=K(1+|x|^2)^{-1}$ with $K>0$ sufficiently small. Then, by Lemma~0 in Allegretto 
\cite{alegr},
the operator $S:E\ri E^*\hookrightarrow E$ defined by $S\varphi =V^+\varphi$ is compact and 
selfadjoint. Next, by Theorem~1 in Allegretto \cite{alegr}, there exist infinitely many 
eigenvalues $\lambda_1(S)>\lambda_2(S)\geq\ldots
\geq \lambda_n(S)\geq \ldots\geq 0$ of $S$ with $\lambda_n(S)\ri 0$ as $n\ri\infty$. So, if
$\mu_n$ ($n\geq 1$) are real numbers such that $\mu_i\geq \mu_j$ if $i<j$
and $\sum_{n=1}^\infty |\mu_n|<\infty$ then, by Theorem \ref{th1}, the mapping
$S\longmapsto\sum_{n=1}^\infty
\mu_n\lambda_n(S)$
is convex and lower semicontinuous.

\subsection{A more general framework}
Consider the class $\kdoi$ of linear selfadjoint operators $S:H\ri H$ having a countable family of 
eigenvalues and such that $S$ can be approximated by operators of finite rank. For any operator 
$S\in\kdoi$, passing eventually at a 
rearrangement, let $\lambda_1(S)\geq \lambda_2(S)\geq\ldots \geq\lambda_n(S)\geq\ldots$ denote the 
eigenvalues of $S$.

Fix a family $\mu=(\mu_1,\mu_2,\ldots ,\mu_n,\ldots)$ of real numbers such that $\mu_i\geq\mu_j$ 
if $i<j$. 
Consider the class $\kdoimu$ of operators $S\in\kdoi$ such that the series
$\sum_{n=1}^\infty\mu_n\lambda_n(S)$ converges.

Under these hypotheses, we establish the following infinite dimensional version of the Schur 
convexity property.

\begin{theorem}\label{th2}
The functional $\psi :\kdoimu\ri\RR$ defined by $\psi (S)=\sum_{n=1}^\infty\mu_n\lambda_n(S)$ is 
convex
and lower semicontinuous.
\end{theorem}

\proof By the definition of $\kdoimu$, for any operator belonging to this class there exists
a sequence $(S_n)_{n\geq 1}$ of operators of finite rank such that $\| S_n-S\|\ri 0$ as 
$n\ri\infty$.
So, by Theorem 4.2 in Gohberg-Krein \cite{gohkre}, we have
$\lim_{n\ri\infty}\lambda_j(S_n)=\lambda_j(S)$, for all positive integer $j$. Define, for all 
$m\geq 1$,
$$\psi_m :\kdoimu\ri\RR\,,\qquad\qquad \psi_m (S)=\sum_{j=1}^m\mu_j\lambda_j(S)\,.$$
Therefore
\begin{equation}\label{bre1}
\lim_{n\ri\infty}\sum_{j=1}^m\mu_j\lambda_j(S_n)=\sum_{j=1}^m\mu_j\lambda_j(S)=\psi_m(S)\,.
\end{equation}
On the other hand, since $S\in\kdoimu$,
\begin{equation}\label{bre2}
\lim_{m\ri\infty}\psi_m(S)=\sum_{j=1}^\infty\mu_j\lambda_j(S)=\psi(S)\,.
\end{equation}

Let $S,T\in\kdoimu$ and assume that $S_n,T_n$ are operators of finite rank such that $\|S_n-S\|\ri 
0$ and
$\|T_n-T\|\ri 0$ as $n\ri\infty$. Applying the Schur convexity property we obtain
$$\psi_m(S_n+T_n)\leq\psi_m(S_n)+\psi_m(T_n)\,,\qquad\mbox{for all $m,n\geq 1$}\,.$$
Taking $n\ri\infty$ and using \eq{bre1} we find
$$\psi_m(S+T)\leq\psi_m(S)+\psi_m(T)\,,\qquad\mbox{for all $m\geq 1$}\,.$$
Next, by \eq{bre2}, we deduce that 
$$\psi(S+T)\leq\psi(S)+\psi(T)\,,\qquad\mbox{for all $S,T\in\kdoimu$}\,.$$
A similar argument shows that $\psi$ is positive homogeneous.

The lower semicontinuity of $\psi$ follows with the same arguments as in the proof of Theorem 
\ref{th1}. \qed

\medskip
{\bf Examples.} 1. {\it Schr\"odinger operators with arbitrary potential}. 
Let $H_0$ denote the differential operator $d^2/dx^2$ on $L^2(0,1)$ with the boundary conditions
$u(0)=u(1)=0$ and assume that $V\in L^\infty (0,1)$ is an arbitrary potential. Let $\lambda_n(S)$ 
be the $n$th 
eigenvalue of the operator $S=H_0+V$. Then, by Theorem XIII.82.5 in Reed and Simon \cite{reed4},
\begin{equation}\label{erwin4}
\lambda_n(S)=-n^2\pi^2+\int_0^1V(x)dx+o(1)\qquad\mbox{as $n\ri\infty$}.\end{equation}

Fix the real numbers $\mu_n$ ($n\geq 1$) such that $\mu_i\geq \mu_j$ if $i<j$ and the series 
$\sum_{n=1}^\infty\mu_n\lambda_n(S)$ converges. 
Using the asymptotic estimate \eq{erwin4}, we deduce that, for the last purpose, it is enough to 
choose $\mu_n$ so that $\mu_n=O\left(n^{-p}\right)$, for some $p>3$. 
Then, by Theorem \ref{th2}, the mapping
$$S\longmapsto\sum_{n=1}^\infty\mu_n\lambda_n(S)$$
is convex and lower semicontinuous.

\smallskip
2. {\it Wave functions on infinite depth wells}. Fix arbitrarily the positive numbers $a$ and $b$. 
Define the
following discontinuous potential energy of a particle in the force field
$$V(x)=\left\{
\begin{array}{lll}
&\di -\infty\qquad &\mbox{if $x<-b$}\\
&\di 0\qquad &\mbox{if $-b<x<a$}\\
&\di -\infty\qquad &\mbox{if $x>a$}\,.
\end{array}\right.$$
Consider the Schr\"odinger equation 
$$\left\{
\begin{array}{ll}
&\di \frac{\hbar^2}{2m}\psi ''+V(x)\psi =\lambda\psi \\
&\di \psi (-b)=\psi (a)=0\,,
\end{array}\right.$$
where $m$ is the mass of the particle and $\hbar$ is Dirac's constant (reduced Planck's constant). 
Cf. Pluvinger \cite[p.~102]{pluvi}, the definition
of $V$ forces $\psi=0$ outside $(-b,a)$.
A straightforward computation shows that the eigenvalues of the associated operator $S$ are given 
by
$$\lambda_n(S)=-\frac{\hbar^2\pi^2}{2m(a+b)^2}\, n^2\,.$$

Fix the real numbers $\mu_n$ ($n\geq 1$) such that $\mu_i\geq \mu_j$ if $i<j$ and the series 
$\sum_{n=1}^\infty\mu_n\lambda_n(S)$ converges. 
The above expression of eigenvalues shows that it is enough to choose $\mu_n$ so that 
$\mu_n=O\left(n^{-p}\right)$, for some $p>3$. 
Applying Theorem \ref{th2}, we deduce that the mapping
$$S\longmapsto\sum_{n=1}^\infty\mu_n\lambda_n(S)$$
is convex and lower semicontinuous.

\smallskip
3. {\it Linear harmonic oscillator}. Consider the Schr\"odinger equation on the whole real axis
\begin{equation}\label{dirac7}\left\{
\begin{array}{ll}
&\di \frac{\hbar^2}{2m}\psi ''+V(x)\psi =\lambda\psi \\
&\di \lim_{|x|\ri\infty}\psi (x)=0=0\,.
\end{array}\right.\end{equation}
In the particular case where $V(x)=-m\omega^2x^2/2$ the above problem describes the linear 
harmonic oscillator. Cf. 
Pluvinger \cite[p.~74]{pluvi} the energy levels of the corresponding linear operator $S$ are given 
by
$\lambda_n(S)=-\hbar\omega(n+1/2)$. So, letting $(\mu_n)_{geq 1}$ so that $\mu_i\geq \mu_j$ if 
$i<j$ and such that the series $\sum_{n=1}^\infty\mu_n\lambda_n(S)$ converges, Theorem \ref{th2} 
implies that the mapping
$S\longmapsto\sum_{n=1}^\infty\mu_n\lambda_n(S)$
is convex and lower semicontinuous. 

We point out that in the case of Morse potentials $V(x)=V_0\left(e^{-2x/a}-2e^{-x/a}\right)$ the 
number of eigenvalues of the
problem \eq{dirac7} is finite.

\smallskip
4. {\it Periodic standing waves of Schr\"odinger's equation}.
In his Ph.D. thesis defended in 1923, de Broglie showed that an electron, or any other particle, 
has a wave associated with it. The second equation established by de Broglie establishes that the 
kinetic energy of a particle is directly proportional to its angular frequency. De Broglie's work 
resulted in the equation $\lambda=\hbar\omega$, where $\lambda$ is the kinetic energy of the 
associated wave and $\omega$ is the angular frequency of the particle. With the same notations as 
in the previous example, we consider the Schr\"odinger equation with periodic boundary conditions
$$\left\{
\begin{array}{ll}
&\di \frac{\hbar^2}{2m}\psi ''+V(x)\psi =\lambda\psi \qquad\mbox{in $(-b,a)$}\\
&\di \psi (-b)=\psi (a)\\ 
&\di \psi' (-b)=\psi ' (a) \,.
\end{array}\right.$$
Outside the fundamental segment of length $L=a+b$ the standing wave $\psi$ is prolonged by 
periodicity such that
$\psi (x+L)=\psi (x)$, for all $x\in\RR$. In Pluvinger \cite[p.~108]{pluvi}
it is provided a class of potentials $V$ for which the associated bound state energies to the 
above problem are given by 
$$\lambda_n(S)=-\frac{2\hbar\pi}{L}\, n \,.$$
Thus, by Theorem \ref{th2}, the mapping $S\longmapsto\sum_{n=1}^\infty\mu_n\lambda_n(S)$
is convex and lower semicontinuous, provided
$(\mu_n)_{n\geq 1}$ are chosen so that $\mu_i\geq \mu_j$ if $i<j$ and the series 
$\sum_{n=1}^\infty\mu_n\lambda_n(S)$ converges.

\smallskip
5. {\it Generalized model of the helium atom}. Let $S$ be the differential operator on 
$L^2(\RR^{3n})$ given by
$$S=\sum_{i=1}^{3n}\left(-\frac{\Delta_i}{2m_i}-\frac{n}{m_i}\right)+\sum_{i<j}\left(\frac{\nabla_
i\cdot\nabla_j}{M}+\frac{1}{|r_i-r_j|}\right)\,,$$
where $M$ and $m_i$ ($1\leq i\leq n$) are arbitrary positive numbers. Cf. Reed and Simon 
\cite{reed4}, the above operator has been introduced by Zhislin and $S$ can be viewed as the 
Hamiltonian of a system consisting of a
nucleus of mass $M$ and $n$ electrons of masses $m_1,\ldots ,m_n$, after the center of the mass 
motion has been removed.
This model generalizes the elementary model of the helium atom which is described by the operator 
$S$ on
$L^2(\RR^6)$ given by
$$S=-\Delta_1-\Delta_2-\frac{2}{|r_1|}-\frac{2}{|r_2|}+\frac{1}{|r_1-r_2|}\,.$$
In both cases (see Kato's Theorem and Theorem XIII.7 in Reed and Simon \cite[p.~89]{reed4}) the 
operator $S$ has a countable family of eigenvalues which can be supposed to be arranged so that 
$\lambda_i(S)\geq\lambda_j(S)$ if $i<j$
(notice that $\lambda_1(S)<-1$ in the case of the elementary model of the helium atom).
Fix the real numbers $\mu_n$ ($n\geq 1$) such that $\mu_i\geq \mu_j$ if $i<j$ and the series 
$\sum_{n=1}^\infty\mu_n\lambda_n(S)$ converges. 
Thus, by Theorem \ref{th2}, the mapping
$S\longmapsto\sum_{n=1}^\infty\mu_n\lambda_n(S)$
is convex and lower semicontinuous. 

\smallskip
6. {\it Schr\"odinger operators with unbounded potential}. Let $V\in L^1_{\rm loc}(\RR^N)$ 
belonging to the class
of operators which are bounded from above and such that $V(x)\ri -\infty$ as $|x|\ri\infty$. Then, 
by Theorem XIII.67
in Reed and Simon \cite{reed4}, the Schr\"odinger operator $S=-\Delta +V$ has a countable family 
of eigenvalues such that
$$\lambda_1(S)\geq\ldots\geq\lambda_n(S)\geq\ldots\qquad\mbox{and $\lambda_n(S)\ri-\infty$ as 
$n\ri\infty$}\,.$$
Consider the real numbers $\mu_n$ ($n\geq 1$) such that $\mu_i\geq \mu_j$ if $i<j$ and the series 
$\sum_{n=1}^\infty\mu_n\lambda_n(S)$ converges. 
Applying Theorem \ref{th2}, we deduce that the mapping
$S\longmapsto\sum_{n=1}^\infty\mu_n\lambda_n(S)$
is convex and lower semicontinuous.

\smallskip
7. {\it Quasilinear anisotropic Sturm-Liouville problems}. Let $\alpha\geq 0$, $p>1$, and $0\leq 
a<b<\infty$. Assume that
$q, s\in L^\infty (a,b)$ and ${\rm ess}\,\inf_{x\in (a,b)}s(x)>0$. Consider the quasilinear 
anisotropic eigenvalue problem
\begin{equation}\label{sturmli}
\left\{\begin{array}{ll}
&\di r^{-\alpha}\left(r^\alpha |u'|^{p-2}u'\right)'+q(r)|u|^{p-2}u=\lambda s(r)|u|^{p-2}u\qquad
\di\mbox{in $(a,b)$}\\
&\di \gamma_1\left(|u|^{p-2}u\right)(a)+ \gamma_2\left(r^\alpha |u'|^{p-2}u'\right)(a)=0\\
&\di \gamma_3\left(|u|^{p-2}u\right)(b)+ \gamma_4\left(r^\alpha |u'|^{p-2}u'\right)(b)=0\,,
\end{array}\right.
\end{equation}
where $\gamma_i\in\RR$ ($i=1,\ldots ,4$) such that $\gamma_1^2+\gamma_2^2>0$ and 
$\gamma_3^2+\gamma_4^2>0$. 

We distinguish two cases: the regular case where $a>0$ or $a=0$ and $0\leq\alpha <p-1$, and the 
singular case defined by $a=0$, $\alpha\geq p-1$. In the singular case the boundary condition at 
the origin is $u'(0)=0$. In both cases Walter
\cite{walter} proved that problem \eq{sturmli} has a countable number of simple eigenvalues 
$\lambda_1(S)>\ldots >
\lambda_n(S)>\ldots$, $\lim_{n\ri\infty}\lambda_n(S)=-\infty$ and the corresponding eigenfunction 
$u_n$ has $n-1$ simple zeroes in $(a,b)$. Consider the real numbers $\mu_n$ ($n\geq 1$) such that 
$\mu_i\geq \mu_j$ if $i<j$ and the series $\sum_{n=1}^\infty\mu_n\lambda_n(S)$ converges. 
So, by Theorem \ref{th2}, the mapping
$S\longmapsto\sum_{n=1}^\infty\mu_n\lambda_n(S)$
is convex and lower semicontinuous.

\medskip
{\bf Conclusions.} In this paper we have extended the Schur convexity property of the eigenvalues
of a symmetric matrix with real entries in the framework of infinite dimensional Hilbert spaces.
First, we have considered the case of linear, selfadjoint, and compact operators. Next, we have
established a corresponding version of the Schur convexity property for linear selfadjoint 
operators that can be approximated by operators of finite rank
and having a countable family of eigenvalues. Our abstract results have been illustrated by 
various examples, including Sturm-Liouville problems, Schr\"odinger operators with variable 
potential, the electron atom model, the linear harmonic oscillator, the generalized model of the 
helium atom, and 
wave functions on infinite depth wells. We have been concerned with linear operators with 
discrete spectrum and our results do not cover the case of operators with a continuous spectrum.

\medskip
{\bf Acknowledgments.} This paper has been written while V. R\u adulescu was visiting
the Laboratoire de M\'ecanique des Solides, Universit\'e de Poitiers, in June 2006. 
V.~R\u adulescu was also partially supported by grants CNCSIS 308/2006 and GAR 80/2006.

\end{document}